    \documentclass[english]{article}
    \newcommand{\mytitle}{Weyl Groups with Coxeter Presentation and
      Presentation by Conjugation} 
    \usepackage{amssymb,amsmath,latexsym}
    \usepackage{eucal}
\renewcommand{\emptyset}{\varnothing}
\newcommand{\newop}[2]{\newcommand{#1}{\mathop{\mathsf{\strut #2}}\nolimits}}
\newop{\Hom}{\mathrm{Hom}}
\newop{\Aut}{\mathrm{Aut}}
\newop{\End}{\mathrm{End}}
\newop{\co}{\mathrm{co\,}}
\newop{\id}{\mathrm{id}}
\newop{\im}{\mathrm{im\,}}
\newop{\rk}{\mathrm{rk}}
\newop{\spann}{\mathrm{span}}
\newop{\conv}{\mathrm{conv}}
\newop{\rank}{\mathrm{rank}}
\newop{\GL}{\mathrm{GL}}
\newop{\PDS}{\mathrm{PDS}}
\newop{\Oo}{\mathrm{O}}
\newop{\cone}{\mathrm{cone}}
\newop{\algint}{\mathrm{algint}}
\newop{\inn}{\mathrm{int}}
\newop{\aff}{\mathrm{aff}}

\newcommand{\R}{\mathbb{R}}

\newcommand{\Z}{\mathbb{Z}}
\newop{\M}{\mathbf{M}}

\newop{\PGL}{\mathrm{PGL}}
\newop{\sL}{\mathrm{sl}}
\newop{\tr}{\mathrm{tr}}

\usepackage{theorem}
\newtheorem{lem}{Lemma}[section]
\newtheorem{thm}[lem]{Theorem}
\newtheorem{prop}[lem]{Proposition}

{\theorembodyfont{\rm}
  \newtheorem{defn}[lem]{Definition}
  \newtheorem{rem}[lem]{Remark}
  \newtheorem{exmp}[lem]{Example}

  \newtheorem{constr}[lem]{Construction}
  }
\newenvironment{proof}{\par\medskip\noindent\upshape\rmfamily
  \textbf{Proof}.\ }{\hspace*{\fill}\rule{1.2ex}{1.2ex}\par\medskip}
\newenvironment{proof*}{\par\medskip\noindent\upshape\rmfamily
  \textbf{Proof}.\ }{\par\medskip}
\newenvironment{lastequation*}{
  \par\medskip\noindent\hspace*{\fill}\begin{math}\displaystyle
  }{
  \end{math}\hspace*{\fill}\rule{1.2ex}{1.2ex}\par}
\newenvironment{lasteqnarray*}{
  \par
  \begin{minipage}[t]{\hfill}
    \begin{eqnarray*}
    }{      
    \end{eqnarray*}
  \end{minipage}
  \begin{minipage}[b]{0.5cm}
    \rule{1.2ex}{1.2ex}
  \end{minipage}
}

    \renewcommand{\tilde}{\widetilde}
    \renewcommand{\hat}{\widehat}
    \renewcommand{\phi}{\varphi}
    \newcommand{\sh}{\mathrm{sh}}
    \newcommand{\ord}{\mathrm{ord}}
    \newcommand{\IRC}{\mathrm{IRC}}
    \newcommand{\trim}{\mathrm{trim}}
    \newcommand{\lng}{\mathrm{lg}}
    \newcommand{\ex}{\mathrm{ex}}
    \newcommand{\calV}{\mathcal{V}}
    \newcommand{\calW}{\mathcal{W}}
    \newcommand{\calS}{\mathcal{S}}
    \begin{document}

    \title{\mytitle}
    \author{Georg Hofmann}
    
    \maketitle
    
    \pagestyle{headings}
    
    \begin{abstract}
      We investigate which Weyl groups have a Coxeter presentation and
      which of them at least have the presentation by conjugation with
      respect to
      their root system. For most concepts of root systems the Weyl group
      has both. In the context of extended affine root systems (EARS) 
      there is a small subclass allowing a Coxeter presentation of
      the Weyl group and a larger subclass allowing the presentation by
      conjugation. We give necessary and sufficient conditions for both
      classes. Our results entail that every extended affine Weyl group
      (EAWeG) 
      has the presentation by conjugation with respect to a suitable EARS.
    \end{abstract}

    \section{Introduction}
    
    If $\calW$ 
    is a group acting on a set $R$ then $R$ is called a 
    \emph{root set} if
    the following conditions are satisfied:
    \begin{enumerate}
\item For every $\alpha\in R$ there is an involution 
      $r_\alpha\in\calW$, and the elements $(r_\alpha)_{\alpha\in R}$
      generate $\calW$.
\item We have 
      \begin{equation*}
        r_\alpha r_\beta(r_\alpha)^{-1}=r_{r_\alpha(\beta)}
      \end{equation*}
      for all $\alpha$ and $\beta\in R$.
    \end{enumerate}        
    The elements of $R$ are called 
    \emph{roots}.
    \index{roots}%
    For any root $\alpha$ we call $r_\alpha$ the involution
    \emph{associated}\index{associated}
    to $\alpha$. The group $\calW$ is referred to as the 
    \emph{Weyl group}
    of $R$
    and group elements of the form $r_\alpha$ are called
    \emph{root involutions}.

    We have introduced the new notion of root sets 
    to encompasses different
    concepts of root systems. For instance $R$ is a root set for the
    Weyl group involved if $R$ 
    \begin{enumerate}
\item[a)] 
      is a finite (not necessarily
      crystallographic) root system. 
\item[b)] 
      consists of the roots associated to
      a \emph{set of root data} in the sense of {\sc Moody} and
      {\sc Pianzola} (see \cite{MoP}), or if $R$
      is the root system of a Coxeter group in the sense of 
      {\sc Deodhar} (see \cite{coxRootSys}) 
      or a root system of the slightly more 
      general kind studied in \cite{diplomarbeit} based on
      {\sc Vinberg}'s 
      idea of a 
      \emph{discrete linear group $\calW$ generated by reflections}
      (see \cite{Vin}).
\item[c)] 
      consists of the real roots of a locally finite root
      system in the sense of \cite{locFinRootSys}.
\item[d)]
      consists of the anisotropic roots of an extended affine
      root system (EARS).
    \end{enumerate}
    We will say that $\calW$ has a 
    \emph{Coxeter presentation with respect to the root set}, 
    if there is a
    subset $\calS$ of
    root involutions in $\calW$ such that $(\calW,\calS)$ is a
    Coxeter system. 

    In cases a) and b) the Weyl group has a Coxeter presentation with
    respect to the root set.
    (See \cite{MoP}~Proposition~5.7.(i), 
    \cite{coxRootSys}~Proposition~3.1.(i) and
    \cite{Vin}~$\S3$~Theorem~2.1 for respective results for the 
    case~b).) We will show for the case d) that the extended affine
    Weyl group (EAWeG) $\calW$ has a Coxeter presentation with respect
    to the root set if and only
    if the nullity of the EARS is less or equal to 1.   

    We will investigate whether the Weyl group 
    $\calW$ of a root set $R$ has
    the presentation 
    \begin{eqnarray}\label{eq:presByConj}
      \big\langle(\hat r_\alpha)_{\alpha\in R^\times}
      &|&
      \hat r_\alpha=\hat r_\beta
      \text{~~~if $\alpha$ and $\beta$ are linearly dependent,}\\
      &&\hat r_\alpha^2=1,\nonumber\\
      &&\hat r_\alpha\hat r_\beta\hat r_\alpha^{-1}
      =\hat r_{r_\alpha(\beta)}; 
      \text{ for }\alpha,\beta\in R^\times\big\rangle,\nonumber
    \end{eqnarray} 
    which is called the 
    \emph{presentation by conjuagtion}.
    The answer is yes in the cases a) and b) since the existence of a
    Coxeter presentation with respect to the root set 
    implies the existence of
    the presentation by conjugation. The answer is also yes in the
    case of c) (see \cite{locFinRootSys}~Theorem~5.12). 

    The answer is also known to be yes for various subclasses of EARSs
    and their EAWeGs (see~\cite{EAWG}, 
    \cite{azamReduced}, \cite{krylyuk}).
    We show that an EAWeG has a presentation by conjugation
    with respect to its EARS
    if and only if it is \emph{minimal} (see
    Definition~\ref{defn:minimalRootSet}.) As a consequence we get the
    following result: Every EAWeG has the presentation by conjugation
    with respect to a suitable EARS.
    
    It has come to our attention that a similar result has been
    obtained independently in \cite{azamA1}. There it is proved for
    the sub-case of EARSs of type $A_1$ that
    the minimality of a certain generating set of $\calW$ is equivalent
    to the existence of the presentation by conjugation. 

    We want to point out that the presentation in (\ref{eq:presByConj})
    is equivalent to the
    presentation 
    \begin{eqnarray*}
      \big\langle(r_\alpha)_{\alpha\in R^\times}
      &|&
      r_\alpha=r_\beta
      \text{~~~if $\alpha$ and $\beta$ are linearly dependent,}\\
      &&r_\alpha r_\beta r_\alpha
      =r_{r_\alpha(\beta)}; 
      \text{ for }\alpha,\beta\in R^\times\big\rangle.
    \end{eqnarray*}

    The author would like to thank Erhard Neher at the University of
    Ottawa for asking the right questions to stimulate the research
    for this work.

    \section[The Coxeter Presentation]%
            {The Coxeter Presentation}
    
    In this section we introduce the definition of a root set. This
    new notion is designed to encompass different concepts of root
    systems as well as the concept of a Coxeter system. The main
    result in this section is that an EAWeG as a Coxeter presentation
    with respect to its EARS if and only if the nullity is greater
    than 1. We postpone the definition of an EARS until Section~3.
            
    \begin{defn}
      Let $\calW$ be a group acting on a set $R$. We call $R$ a 
      \emph{root set}\index{root set}
      for $\calW$ if the following conditions are satisfied
      \begin{enumerate}
  \item For every $\alpha\in R$ there is an involution 
        $r_\alpha\in\calW$, and the elements $(r_\alpha)_{\alpha\in R}$
        generate $\calW$.
  \item We have 
        \begin{equation*}
          r_\alpha r_\beta(r_\alpha)^{-1}=r_{r_\alpha(\beta)}
        \end{equation*}
        for all $\alpha$ and $\beta\in R$.
      \end{enumerate}        
      If $R$ is a root set for $\calW$, then we define
      \begin{equation*}
        \alpha\sim\beta:\iff r_\alpha=r_\beta.
      \end{equation*}
      The elements of $R$ are called 
      \emph{roots}.
      \index{roots}%
      For any root $\alpha$ we call $r_\alpha$ the reflection
      \emph{associated}\index{associated}
      to $\alpha$. The group $\calW$ is referred to as the 
      \emph{Weyl group}
      of $R$
      and group elements of the form $r_\alpha$ are called
      \emph{root involutions}.
    \end{defn}

    \begin{rem}
      The relation $\sim$ is a congruence relation, i.e. an
      equivalence relation satisfying
      \begin{equation*}
        \alpha\sim\beta\iff w.\alpha\sim w.\beta,
      \end{equation*}
      for every $w\in\calW$.
    \end{rem}

    \begin{exmp}\label{exmp:rootSets}
      We refer the reader to the list of examples a)~-~d) given in the
      introduction. We add one example to this list:
      \begin{enumerate}
  \item[e)] If $(\calW,\calS)$ is a Coxeter system, then 
        \begin{equation}\label{eq:t}
          T:=\{wsw^{-1}:w\in W,s\in\calS\}
        \end{equation}
        is a root set for $\calW$ if we set $r_t=t$. The equivalence
        relation in this setting is equality. 
      \end{enumerate}
    \end{exmp}

    Let $R$ be a root set with Weyl group $\calW$. Let $S$ be a subset of
    $R$ such that the group root involutions associated to the
    elements in $S$ generate $W$ and such that
    \begin{equation*}
      (\forall\alpha,\beta\in S)~\alpha\sim\beta\implies\alpha=\beta.
    \end{equation*}
    Set
    \begin{eqnarray}\label{eq:coxPres}
      \hat\calW_S:=
      \big\langle(\hat r_\alpha)_{\alpha\in S}
      &|&
      \hat r_\alpha^2=1        
      \text{~~~for $\alpha\in S$,}\\
      &&(\hat r_\alpha\hat r_\beta)^{\ord(r_\alpha r_\beta)}=1
      \text{~~~for $\alpha$, $\beta\in R$ 
        and $\ord(r_\alpha r_\beta)<\infty$}
      \big\rangle\nonumber
    \end{eqnarray}
    where $\ord(\cdot)$ stands for the order of a group element.
    Since the relations in this presentation are true for the
    generators $(r_\alpha)_{\alpha\in S}$ of $\calW$, there is a
    group homomorphism $\phi:\hat\calW_S\to\calW$ 
    that maps $\hat r_\alpha$ to $r_\alpha$ for all $\alpha\in S$. 
    Since $\phi$ has a generating set of $\calW$ in its image, it is
    surjective. 

    \begin{defn}
      We say that $\calW$ has a
      \emph{Coxeter presentation (with respect to $R$ and $S$)}%
      \index{Coxeter presentation}
      if $\phi$ is injective. 
    \end{defn}

    The requirement in this definition is equivalent to saying that
    there is a subset $\calS$ of root involutions such that 
    $(\calW,\calS)$ is a
    Coxeter system.

    \begin{exmp}
      In the list of examples discussed in the Introduction and in
      Example~\ref{exmp:rootSets} 
      the following Weyl groups have a Coxeter
      presentation with respect to their root sets: a), b), e).
    \end{exmp}

    We will present an example of a root set $R$ with a Weyl group
    which does not permit a Coxeter presentation with respect to $R$.  

    \begin{exmp}\label{exmp:nonCox}
      Consider the real vector space $V:=\R^5$ with the bilinear form
      $(\cdot,\cdot)$ with matrix
      \begin{equation*}
        \begin{pmatrix}
          0&0&0&1&0\\
          0&0&0&0&1\\
          0&0&1&0&0\\
          1&0&0&0&0\\
          0&1&0&0&0
        \end{pmatrix}
      \end{equation*}
      with respect to the standard basis. 
      In $V$ we consider the set
      \begin{eqnarray*}
        R^\times
        &:=&
        \big\{(z_1,z_2,z_3,0,0)|z_1,z_2\in\Z,~z_3\in\{1,-1\},
        ~z_1z_2\text{ is even}\big\}.
      \end{eqnarray*}
      (This is the set of anisotropic roots of an EARS 
      of type $A_1$ with nullity 2.) 

      To an element $\alpha\in R^\times$ we associate the reflection
      \begin{equation*}
        r_\alpha:~V\to V,~v\mapsto v-2\frac{(v,\alpha)}{(\alpha,\alpha)}\alpha.
      \end{equation*}
      Let $\calW$ be the subgroup of
      $\Aut(V)$ generated by the reflections 
      $\{r_\alpha,\alpha\in R^\times\}$. Then $R^\times$ is a root set
      for $\calW$. 
    \end{exmp}    
    
    \begin{lem}
    The group $\calW$ in the previous example does not have a 
    Coxeter presentation with respect to
    $R^\times$. 
    \end{lem}
    
    \begin{proof}
      We will follow the proof in \cite{myThesisBook}~Example~4.3.36.
      Suppose that $S$ is a subset of $R^\times$ such that $\calW$ has a
      Coxeter presentation with respect to $R^\times$ and $S$. 
      Then  $\spann S$ is three dimensional: The subspace 
      $U:=\spann S$ has at least dimension 1, since $S$ contains an
      element of $R^\times$. Now $U$ is a $\calW$-invariant
      subspace. Dimension 1 or 2 contradict 
      the fact the orbit under $\calW$ of any root spans
      a 3-dimensional subspace.
      
      It is easily verified that the product of two reflections in 
      $\calW$ has infinite order. 
      So $\hat\calW_S$ is the free group generated by
      the involutions associated to $S$. We showed that $S$ must contain 3
      linearly independent roots say $\alpha_1$, $\alpha_2$ and
      $\alpha_3$. Denote the subgroup  
      of $\calW$
      generated by the set 
      $\calS':=\{r_{\alpha_1},r_{\alpha_2},r_{\alpha_3}\}$ by
      $\calW'$. 
      We are done if we show
      that $(\calW',\calS')$ is not a Coxeter system. 
      (For the
      necessary result about parabolic subgroups see
      \cite{Hum}~Theorem in Section~5.5).  
      So 
      without loss of generality, we may assume that we have
      \begin{equation*}
        \alpha_1=(0,0,1,0,0)^T,~
        \alpha_2=(1,0,1,0,0)^T,~\text{and}~
        \alpha_3=(0,1,1,0,0)^T.
      \end{equation*}
      
      For the associated reflections $r_i:=r_{\alpha_i}$ ($i=1$, 2, 3)
      it can be verified by geometric arguments or 
      by a matrix computation that
      \begin{eqnarray*}
        r_1r_2r_3r_1r_2r_3r_2r_1r_3r_2r_1&=&r_3
        \text{~~~or, equivalently}\\
        r_1r_2r_3r_1r_2r_3r_2r_1r_3r_2r_1r_3&=&1
      \end{eqnarray*}
      We provide here the matrices corresponding to the three reflections
      $r_1$, $r_2$, $r_3$ and leave the computations to the reader:
      \begin{equation*}
        \begin{pmatrix}
          1& 0& 0& 0& 0\\
          0& 1& 0& 0& 0\\
          0& 0&-1& 0& 0\\
          0& 0& 0& 1& 0\\
          0& 0& 0& 0& 1
        \end{pmatrix}
        ,~~
        \begin{pmatrix}
          1& 0&-2&-2& 0\\
          0& 1& 0& 0& 0\\
          0& 0&-1&-2& 0\\
          0& 0& 0& 1& 0\\
          0& 0& 0& 0& 1
        \end{pmatrix}
        ,~~
        \begin{pmatrix}
          1& 0& 0& 0& 0\\
          0& 1&-2& 0&-2\\
          0& 0&-1& 0&-2\\
          0& 0& 0& 1& 0\\
          0& 0& 0& 0& 1
        \end{pmatrix}
      \end{equation*}
      
      The relation above is not satisfied in $\hat\calW_S$.
      See for example, \cite{Hum}~8.1: In $\hat\calW_S$ 
      any word of
      length greater than zero, which contains no subsequent repetition of
      an element of $\{\hat r_\alpha:\alpha\in S\}$, 
      does not represent the identity.
      So the homomorphism $\hat\calW_S\to\calW$ is not injective.
    \end{proof}
    
    In \cite{myThesisBook}~Lemma~4.3.33 it is shown that the involutions
    in $\calW$ are exactly the reflections in $\calW$. Moreover, the
    reflections in $\calW$ are exactly the root reflections, 
    in other words there are no ghost
    reflections, i.e. reflections in the group that are not associated to
    a root. 
    (From \cite{myThesisBook}~Theorem~4.4.47 and Theorem~4.4.51 it follows
    that these ghost reflections only occur in higher dimensions).
    
    All of this taken together results in
    \begin{prop}\label{prop:NonCox}
      The group $\calW$ is not a Coxeter group, i.e. there is no root set
      $R$ for $\calW$ such that $\calW$ has a Coxeter presentation with
      respect to $R$.
    \end{prop}

    The question arrises which EARS have a Coxeter presentation of the
    Weyl group with respect to the underlying root set.    
    An EARS with nullity 0 is just a finite root
    system with zero added. It is known that a Coxeter presentation
    with respect to the root set
    exists for this case. The anisotropic roots of an EARS with
    nullity 1 can be viewed as a root system in the sense of
    \cite{MoP}. Again, a Coxeter presentation with respect to
    the underlying root set exists. If the
    nullity of an EARS is greater than or equal to 2 than there are always
    three anisotropic roots $\alpha_1$, $\alpha_2$ and $\alpha_3$ 
    whose reflections
    generate a subgroup of the Weyl group isomorphic to $\calW$ in the
    previous example. This means that the Weyl group does not allow a
    Coxeter presentation with respect to the underlying root set, 
    since a subgroup of a Coxeter group generated by reflections is
    always a Coxeter group (See \cite{Deodhar}, \cite{Dyer} or 
    \cite{myThesisBook}~Corollary~2.3.29). So we have

    \begin{thm}
      Let $\calW$ be the Weyl group of an EARS with nullity $\nu$. Denote
      the set of anisotropic roots by $R^\times$. 
      Then $\calW$ has a Coxeter presentation with respect to $R^\times$
      if and only if $\nu<2$.
    \end{thm}
    
    It is more difficult to judge whether a given Weyl group of an EARS of
    nullity $\nu\ge2$ has a Coxeter presentation with respect to any
    root set. The arguments used for the previous
    example to obtain Proposition~\ref{prop:NonCox} 
    may fail in two ways: 
    Firstly, not all involutions
    of the Weyl group must be reflections: The finite root system $B_2$
    already contains a rotation about $180^\circ$. Secondly, some
    reflections may be ghost reflections, i.e. they are not associated to
    a root in the root system.

    \section[The Presentation by Conjugation]%
            {The Presentation by Conjugation}

    In this section we introduce the notion of the presentation by
    conjugation of a group with respect to a root set. If a group
    $\calW$ has a Coxeter presentaion with respect to a root set $R$
    then $\calW$ has a presentaion by conjugation with respect to this
    root set. The key result of this section is that if a $\calW$ has
    the presentation by conjugation with respect to $R$ then $R$ is
    minimal. 

    Let $R$ be a root set with Weyl group $\calW$. Set
    \begin{eqnarray}
      \hat\calW:=
      \big\langle(\hat r_\alpha)_{\alpha\in R}
      &|&
      \hat r_\alpha=\hat r_\beta
      \text{~~~for $\alpha$ and $\beta\in R$ with $\alpha\sim\beta$,}
      \nonumber\\
      &&\hat r_\alpha^2=1        
      \text{~~~for $\alpha\in R$,}
      \label{eq:preConj}\\
      &&\hat r_\alpha \hat r_\beta \hat r_\alpha^{-1}
      =\hat r_{r_\alpha(\beta)}
      \text{~~~for $\alpha$ and $\beta\in R$}
      \big\rangle.\nonumber
    \end{eqnarray}
    Since the relations in this presentation are true for the
    generators $(r_\alpha)_{\alpha\in R}$ of $\calW$, there is a
    group homomorphism $\phi:\hat\calW\to\calW$ 
    that maps $\hat r_\alpha$ to $r_\alpha$ for all $\alpha\in R$. 
    Since $\phi$ has a generating set of $\calW$ in its image, it is
    surjective. 

    \begin{defn}
      We say that $\calW$ has the
      \emph{presentation by conjugation (with respect to $R$)}%
      \index{presentation by conjugation}
      if $\phi$ is injective. 
    \end{defn}

    The following result is known. Since no source 
    for it is known to the author we will include a proof.
    \begin{prop}
      Let $(\calW,S)$ be a Coxeter system. Consider the root set $T$ for
      $W$ described in Example~\ref{exmp:rootSets}. Then $W$ has the
      presentation by conjugation with respect to $T$.
    \end{prop}

    \begin{proof}
      Set $\hat S=\{\hat r_s:s\in S\}$. We want to show that the map
      $S\to\hat S,~s\mapsto\hat r_s$ can be extended to a map 
      $\psi:\hat\calW\to\calW$. For that purpose we need 
      to show that the relations defining the Coxeter group $\calW$
      (see (\ref{eq:coxPres})) 
      hold for the generators 
      $\hat S$ in $\hat\calW$. The first defining relation in
      (\ref{eq:coxPres}) is clear so let's look at the second:

      In $\calW$ we have:
      \begin{eqnarray*}
        (s_1s_2)^{m(s_1,s_2)}&=&1.
      \end{eqnarray*}
      Let's consider the case that $m(s_1,s_2)$ is odd, say
      $m(s_1,s_2)=2k+1)$. Then this relation 
      together with the fact that $s_1$ and $s_2$ are
      involutions entails
      \begin{eqnarray*}      
        (s_1s_2)^ks_1(s_2^{-1}s_1^{-1})^k&=&s_2.
      \end{eqnarray*}
      From the defining relations of $\hat\calW$ it follows that in 
      $\hat\calW$ we have
      \begin{eqnarray*}      
        (\hat r_1\hat r_2)^k\hat r_1(\hat r_2^{-1}\hat r_1^{-1})^k
        &=&
        \hat r_2
      \end{eqnarray*}
      where we have written $\hat r_1$ for $\hat r_{s_1}$ and 
      $\hat r_2$ for $\hat r_{s_2}$ to avaoid an notational overload. 
      In turn, this implies
      \begin{eqnarray*}
        (\hat r_1\hat r_2)^{m(\hat r_1,\hat r_2)}&=&1.
      \end{eqnarray*}      
      The case $m(s_1,s_2)$ odd is similar and is left to the reader.

      So we have the map $\psi:\calW\to\hat\calW$ and need to prove that 
      $\psi\circ\phi$ is the identity on $\hat\calW$. It is clear that it
      is the identity on $\hat S$, so we need to show that $\hat S$
      generates $\hat\calW$. We do this by proving that the subgroup 
      $\langle\hat S\rangle$ of $\hat\calW$ generated by $\hat S$
      contains the set
      $\hat T:=\{\hat r_t:t\in T\}$.

      So let $\hat r_t\in \hat T$. this means $t=wsw^{-1}$ for some
      $w\in\calW$ and $s\in S$. The group element $w$ can be written as 
      $w=s_1s_2\cdots s_k$ for some $s_1$, $s_2,\dots,s_k\in S$. So we
      have
      \begin{eqnarray*}
        t&=&s_1s_2\cdots s_k\cdot 
        s\cdot 
        s_k^{-1}s_{k-1}^{-1}\cdots s_1^{-1}.
      \end{eqnarray*}
      Using arguments similar to those above we deduce
      \begin{eqnarray*}
        \hat r_t
        &=&
        \hat r_1\hat r_2\cdots \hat r_k
        \cdot 
        \hat r_s
        \cdot 
        \hat r_k^{-1}\hat r_{k-1}^{-1}\cdots \hat r_1^{-1}.
      \end{eqnarray*}
      This shows $\hat r_t\in\langle\hat S\rangle$, 
      which concludes the proof.
    \end{proof}

    \begin{exmp}
      \begin{enumerate}
  \item 
        In the list of examples discussed in the Introduction the root
        sets in a) and b) have a Coxeter presentation. By the previous
        proposition, they have a presentation by conjugation.
  \item 
        The underlying root set of a locally finite root system
        described in the introduction in c) has the
        presentation by conjugation (see
        \cite{locFinRootSys}~Theorem~5.12). 
      \end{enumerate}
    \end{exmp}

    \begin{defn}\label{defn:minimalRootSet}
      Let $R$ be a root set with Weyl group $\calW$. 
      We call the root set $R$ 
      \emph{minimal}\index{minimal}
      if there is no root
      $\gamma\in R$ such that reflections associated to the elements
      of the set
      \begin{equation}
        R_\gamma
        :=R\setminus\{\alpha:(\exists w\in\calW) w.\gamma\sim\alpha\}
      \end{equation}
      generate $\calW$.
    \end{defn}

    \begin{thm}
      If the group $\calW$ has a presentation by conjugation with respect
      to $R$, then the root set $R$ is minimal.
    \end{thm}

    \begin{proof}
      Suppose that $R$ is not minimal, i.e. 
      there is a $\gamma\in R$ such that reflections
      associated to the elements of the set 
      \begin{math}
        R_\gamma
      \end{math}
      in (\ref{eq:rGamma})
      generate $\calW$. Since $r_\gamma$ is an element 
      of the group $\calW$
      there must be a relation of the kind
      \begin{equation}\label{eq:rGamma}
        r_\gamma=r_{\delta_1}r_{\delta_2}\cdots r_{\delta_k}
        ~~\text{with $\delta_1,\dots,\delta_k\in R_\gamma$.}
      \end{equation}
      Now let $F(R)$ be the free group generated by the elements
      $r_\alpha$ with $\alpha\in R$. Consider the normal
      subgroup $N$ of $F(R)$ generated by 
      \begin{eqnarray*}\label{eq:generators2}
        X:=&&\{\hat r_\alpha(\hat r_\beta)^{-1}:\alpha,\beta\in R
        ~\text{with $\alpha\sim\beta$}\}\\
        &\cup&
        \{\hat r_\alpha\hat r_\alpha:
        \alpha\in R\}\\
        &\cup&
        \{\hat r_\alpha\hat r_\beta\hat r_\alpha
        (\hat r_{r_\alpha(\beta)})^{-1}:
        \alpha,\beta\in R\}.
      \end{eqnarray*}
      The quotient $F(R)/N$ is the group $\hat\calW$.
      
      Let's consider the map 
      \begin{eqnarray*}
        \psi:&R&\to~\Z_2,\\
        &\alpha&\mapsto~
        \begin{cases}
          1&\text{if } 
          \alpha\in R_\gamma\\
          0&\text{else}
        \end{cases}
      \end{eqnarray*}
      This map extends to a group homomorphism $\tilde\psi:F(R)\to\Z_2$.
      It is easily verified that $X$ is in $\ker\psi$ and thus
      $N\le\ker\psi$. This means
      that the map $\tilde\psi$ factors to a map $\hat\psi:F(R)/N\to\Z_2$.
      
      Recall the relation (\ref{eq:rGamma}) in $\calW$. This relation is
      not true in $\hat\calW=F(R)/N$ since
      \begin{equation*}
        \hat\psi(\hat r_\gamma)=1 
        ~~~\text{whereas}~~~
        \hat\psi(\hat r_{\delta_1}\hat r_{\delta_2}
        \cdots\hat r_{\delta_k})=0.
      \end{equation*}
      This implies that the homomorphism 
      $\phi:\hat\calW\to\calW$ is not injective. 
    \end{proof}

    The converse of the previous theorem is not true. In the following
    we present an example of a reflection group with a minimal root
    system which does not have a presentation by conjugation.

    \begin{exmp}
      Recall the root set $R^\times$ with Weyl group $\calW$ from
      Example~\ref{exmp:nonCox}. Let $U$ be the subspace of $V$
      spanned by the first three standard basis vectors. Then $U$ is
      invariant under $\calW$. Let $\overline\calW$ 
      be the restriction of $\calW$
      to $U$. Since $R^\times\subseteq U$, the set $R^\times$ is a root
      set for $\overline\calW$. 
      The restriction induces a group homomorphism
      \begin{equation*}
        \calW\to\overline\calW.
      \end{equation*}
      (The group $\calW$ is actually a central extension of 
      $\overline\calW$.) This homomorphism is not injective: Set
      \begin{equation*}
        \alpha_0:=
        \begin{pmatrix}
          0\\0\\1\\0\\0
        \end{pmatrix},~
        \alpha_1:=
        \begin{pmatrix}
          2\\0\\1\\0\\0
        \end{pmatrix},~
        \alpha_2:=
        \begin{pmatrix}
          2\\1\\1\\0\\0
        \end{pmatrix},~
        \alpha_3:=
        \begin{pmatrix}
          0\\1\\1\\0\\0
        \end{pmatrix}.
      \end{equation*}
      Then $r_{\alpha_0}r_{\alpha_1}r_{\alpha_2}r_{\alpha_3}$ is a
      nontrivial element in the kernel of the homomorphism. The
      groups $\hat\calW$ and $\hat{\overline\calW}$ have the same
      presentation. 
      Since the homomorphism 
      $\hat\calW\to\calW$ is surjective, the composition of homomorphisms
      \begin{equation*}
        \hat\calW\to\calW\to\overline\calW
      \end{equation*}
      is not injective. This means that 
      $\hat{\overline\calW}\to\overline\calW$ is not injective, so 
      $\overline\calW$
      does not have the
      presentation by conjugation with respect to $R^\times$.

      The orbits of $\overline\calW$ in $R^\times$ are given by
      \begin{equation*}
        \overline\calW.\alpha
        =
        \{\alpha,-\alpha\}
        +\big((2\Z)\times(2\Z)\times\{0\}\big).
      \end{equation*}
      From this, it is not hard to see that $R$ is a minimal root set
      for $\overline\calW$.
    \end{exmp}
    
    \begin{exmp}
      For all concepts of root systems $R$ and their Weyl groups
      $\calW$ such that the underlying root 
      set provides the presentation by
      conjugation for the Weyl group, 
      the root set is minimal. This is true for the examples a)~-~c)
      in the Introduction as well as 
      Example~\ref{exmp:rootSets}~e).
    \end{exmp}
    
    We will look at the case of a finite crystallographic root system
    more closely and derive a result about the orbits of its Weyl
    group which we need later on:
    
    \begin{lem}\label{lem:finRootSub}
      Let $R$ be an irreducible finite crystallographic 
      root system with Weyl group
      $\calW$. Suppose $\tilde R$
      is a subset of $R$ that is invariant under $\calW$ and such that
      the reflections associated to the roots in $\tilde R$ generate
      $\calW$. Then $\tilde R$ is an irreducible finite
      crystallographic
      root system and
      one of the following is true:
      \begin{enumerate}
  \item $R=\tilde R$,
  \item $R$ is of type $BC_1$ and $\tilde R$ is of type $A_1$, 
  \item $R$ is of type $BC_2$ and $\tilde R$ is of type $B_2$, 
  \item $R$ is of type $BC_\ell$, $\ell\ge3$ and $\tilde R$ is of type
        $B_\ell$ or $C_\ell$.
      \end{enumerate}
    \end{lem}
    
    \begin{proof}
      Let $R_\sh$, $R_\lng$ and $R_\ex$ be the set of all short, long
      and extra long roots in $R$. Then $R$ is the disjoint union of
      these three sets. Moreover these three sets are precisely the
      orbits under $\calW$ in the root system. (See
      \cite{bou}~Ch.~VI.~$n^0$~1.4.10.) So $\tilde R$ must be the
      union of some of the three sets.

      Note that with $R_\gamma$ as in 
      Definition~\ref{defn:minimalRootSet} we have
      \begin{equation*}
        R_\gamma=
        \begin{cases}
          R_\lng&\text{if $\gamma\in R_\sh\cup R_\ex$}\\
          R_\sh\cup R_\ex&\text{if $\gamma\in R_\lng$}.
        \end{cases}
      \end{equation*}      
      If $R$ is simply laced, i.e. if $R_\ex$ is empty, then
      $R=\tilde R$, since finite root systems are minimal root
      sets. Otherwise, we have 
      \begin{equation*}
        \tilde R=R\setminus R_\sh
        \text{~~~or~~~}
        \tilde R=R\setminus R_\ex
      \end{equation*}
      for the same reason. 
      Then we are in one of the cases (ii) - (iv). In each of these
      cases $\tilde R$ is an irreducible finite root system.       
    \end{proof}

    In the following we will present another result needed later on to
    prove that certain EAWeGs have the presentation by conjugation.

    Let $R$ be a minimal root system with Weyl group $\calW$. Let 
    $\tilde R$ be a $\calW$-invariant subset 
    of $R$ and denote by $\calW_{\tilde R}$ the
    subgroup of $\calW$ generated by the reflections associated to the
    elements of $\tilde R$. Suppose the following
    property is satisfied: 
    \begin{equation*}
      (\forall\alpha\in R\setminus\tilde R)~
      (\forall w\in\calW)~
      (\exists w'\in\calW_{\tilde R})~
      w\alpha\sim w'\alpha
    \end{equation*}


    \begin{prop}\label{prop:SubRootPres}
      If the group $\calW_{\tilde R}$ has a presentation by conjugation
      with respect to $\tilde R$, 
      then the group $\calW$ has a presentation 
      by conjugation with respect
      to $R$. 
    \end{prop}

    \begin{proof}
      We will start with a relation 
      \begin{equation}
        \label{eq:startRel}
        \hat r_*\hat r_*\cdots \hat r_*
      \end{equation}
      satisfied in the group $\calW$ and
      show that the relation follows from the relations of the
      presentation by conjugation:
      \begin{eqnarray}
        \label{eq:rels1}
        &&
        \hat r_\alpha=\hat r_\beta
        \text{~~~for $\alpha$ 
          and $\beta\in R$ with $\alpha\sim\beta$,}\\
        \label{eq:rels2}
        &&\hat r_\alpha^2=1        
        \text{~~~for $\alpha\in R$,}\\
        \label{eq:rels3}
        &&\hat r_\alpha\hat r_\beta\hat r_\alpha^{-1}
        =\hat r_{r_\alpha(\beta)}
        \text{~~~for $\alpha$ and $\beta\in R$}
      \end{eqnarray}

      First we suppose that relation~(\ref{eq:startRel})
      contains no elements 
     $\hat r_\alpha$ with $\alpha\in R\setminus\tilde R$. 
     But then we know
     by hypothesis that relation~(\ref{eq:startRel})
     follows from the relations 
     in (\ref{eq:rels1}) to (\ref{eq:rels3}) with $\alpha$ and
     $\beta\in\tilde R$. So in this case there is nothing to prove.

     Now we suppose that relation~(\ref{eq:startRel})
     contains at least one $\hat r_\alpha$
     with $\alpha\in R\setminus\tilde R$. Moreover, we assume it
     contains no $\hat r_{\alpha'}$ such that $\alpha\sim w.\alpha'$ for
     some $w\in\calW$. Then the relation can be written as
     \begin{eqnarray*}
       \hat r_\alpha
       &=&
       \hat r_*\hat r_*\cdots\hat r_*,
     \end{eqnarray*}
     where all the generators on the right side are associated to
     roots in 
     \begin{equation*}
       R_\alpha
       =
       R\setminus\{\gamma:(\exists w\in\calW) \gamma\sim w.\alpha\}.
     \end{equation*}
     This means that the reflections associated to this set
     actually generate $\calW$, which is a
     contradiction to the fact that $R$ is minimal.
     
     So, if we are supposing that relation~(\ref{eq:startRel})
     contains an
     $\hat r_\alpha$
     with $\alpha\in R\setminus\tilde R$ then it must contain another
     element $\hat r_{\alpha'}$ with $\alpha\sim w.\alpha'$ for some 
     $w\in\calW$. So it is of the form
     \begin{eqnarray}
       \label{eq:ourRelation}
       \hat r_*\cdots\hat r_*
       \hat r_\alpha
       \hat r_*\cdots\hat r_*
       \hat r_{\alpha'}
       \hat r_*\cdots\hat r_*
       &=&
       1.
     \end{eqnarray}
     
     By hypothesis, we may assume $w\in\calW_{\tilde R}$. So
     we have 
     \begin{eqnarray*}
       w
       &=&
       r_{\alpha_1}
       r_{\alpha_2}
       \cdots
       r_{\alpha_n}
     \end{eqnarray*}
     for some $\alpha_1,\alpha_2,\dots,\alpha_n\in\tilde R$. The
     conjugation relations (\ref{eq:rels1}) and (\ref{eq:rels3}) imply
     \begin{eqnarray*}
       \hat r_\alpha
       &=&
       \hat r_{w.\alpha'}
       =
       \hat r_{
         r_{\alpha_1}
         r_{\alpha_2}
         \cdots
         r_{\alpha_n}
         .\alpha'
       }\\
       &=&       
       \hat r_{\alpha_1}
       \hat r_{\alpha_2}
       \cdots
       \hat r_{\alpha_n}
       ~
       \hat r_{\alpha'}
       ~       
       \hat r_{\alpha_n}^{-1}
       \hat r_{\alpha_{n-1}}^{-1}
       \cdots
       \hat r_{\alpha_1}^{-1},
     \end{eqnarray*}
     which can be rewritten as
     \begin{eqnarray*}
       \hat r_{\alpha'}
       &=&
       \hat r_{\alpha_n}
       \hat r_{\alpha_{n-1}}
       \cdots
       \hat r_{\alpha_1}
       ~\hat r_\alpha~
       \hat r_{\alpha_1}
       \hat r_{\alpha_2}
       \cdots
       \hat r_{\alpha_n}
     \end{eqnarray*}
     using relation (\ref{eq:rels2}).
     So our relation (\ref{eq:ourRelation}) can be written as
     \begin{eqnarray}
       \label{eq:ourRelation2}
       \hat r_*\cdots\hat r_*
       \hat r_\alpha
       \hat r_*\cdots\hat r_*
       \hat r_{\alpha}
       \hat r_*\cdots\hat r_*
       &=&
       1.
     \end{eqnarray}
     Note that number of generators in this relation 
     may have increased but the
     number of generators elements $\hat r_\beta$ with 
     $\beta\in R\setminus\tilde R$ has not changed. 

     For any $\beta\in R$ the conjugation relation 
     $\hat r_\alpha \hat r_\beta \hat r_\alpha^{-1}
     =\hat r_{r_\alpha(\beta)}$
     can be rewritten as 
     \begin{eqnarray*}
       \hat r_\beta\hat r_\alpha
     &=&\hat r_\alpha\hat r_{r_\alpha(\beta)}
     \end{eqnarray*}
     using relation (\ref{eq:rels2}). This means that the relation
     (\ref{eq:ourRelation2}) can be rewritten as
     \begin{eqnarray*}
       \hat r_*\cdots\hat r_*
       \hat r_\alpha
       \hat r_{\alpha}
       \hat r_*\cdots\hat r_*
       &=&
       1.
     \end{eqnarray*}
     The two generators $\hat r_\alpha$ in this relation cancel out
     reducing it to a relation which has fewer generators 
     $\hat r_\beta$ 
     with 
     $\beta\in R\setminus\tilde R$. This process can be repeated until
     the number of generators $\hat r_\beta$ 
     with 
     $\beta\in R\setminus\tilde R$ is zero. Then we are done  since
     we are in a case that has been discussed earlier.
    \end{proof}

    \section{Extended affine root systems}
    
    At the beginning of this section 
    we present the definition of extended affine root
    systems and some basic facts about them taken from
    \cite{EALA_AMS}. The main results of this work are derived after
    that: Given an EARS $R$, the corresponding Weyl group $W$ 
    has the presentation by conjugation with respect to $R$
    if and only if
    $R$ is minimal. Moreover, every EARS $R$ contains a minimal
    one. Taken together, this means that every EAWeG has a
    presentation by conjugation with respect to some EARS.

    \begin{defn}
      Let $\calV$ be a finite dimensional real vector space with a
      positive semi-definite symmetric bilinear form $(\cdot,\cdot)$ and
      let $R$ be a subset of $\calV$. Set
      \begin{equation*}
        \calV^0:=\{v\in\calV:(v,v)=0\},~~
        R^0:=R\cap\calV^0,
        \text{~~and~~}
        R^\times:=R\setminus R^0
      \end{equation*}
      The subset $R$ is called an
      \emph{extended affine root system}
      \index{extended affine!-- root system|uu}%
      in $\calV$ if it satisfies the following axioms:
      \begin{enumerate}
  \item[(R1)]
        $0\in R$.
  \item[(R2)]
        $-R=R$.
  \item[(R3)]
        $R$ spans $\calV$.
  \item[(R4)]
        $\alpha\in R^\times\implies 2\alpha\notin R$.
  \item[(R5)]
        $R$ is discrete in $\calV$.
  \item[(R6)]
        If $\alpha\in R^\times$ and $\beta\in R$, then there exist
        $d,u\in\Z_{\ge 0}$ satisfying
        \begin{equation*}
          \{\beta+n\alpha:n\in\Z\}\cap R
          =
          \{\beta-d\alpha,\dots,\beta+u\alpha\}
          \text{~~and~~}
          d-u=2\frac{(\alpha,\beta)}{(\alpha,\alpha)}.
        \end{equation*}
  \item[(R7)]
        $R^\times$ cannot be decomposed as a disjoint union 
        $R_1\dot\cup R_2$, where $R_1$ and $R_2$ are nonempty subsets of
        $R^\times$ satisfying $(R_1,R_2)=\{0\}$.
  \item[(R8)]
        If $\sigma\in R^0$, 
        then there exists $\alpha\in R^\times$ such that
        $\alpha+\sigma\in R$.
      \end{enumerate}
      Elements of $R^\times$ are called 
      \emph{anisotropic roots}
      \index{anisotropic!-- roots|uu}%
      and elements of $R^0$ are called 
      \emph{isotropic roots}. 
      \index{isotropic!-- root|uu}%
      The \emph{nullity}%
      \index{nullity|uu}
      of $R$ is the dimension of $\calV^0$.
    \end{defn}
    
    By $\overline\cdot:\calV\to\calV/\calV^0$ we denote the quotient map.
    Then $\overline R^\times$ is a finite crystallographic
    root system. It contains a
    fundamental system 
    $\overline\Pi=\{\overline\alpha_1,\dots,\overline\alpha_l\}$ with 
    ${l=\dim(\calV/\calV^0)}$. For each $\overline\alpha_i$ we pick a
    preimage $\dot\alpha_i$ in $R^\times$. Then we
    set
    \begin{equation}\label{eq:dotCalV}
      \dot\calV:=
      \spann_\R\{\dot\alpha_1,\dots,\dot\alpha_l\}.
    \end{equation}
    
    \begin{defn}
      A nonempty subset $S$ of a real vector space $V$ is called a 
      \emph{translated semilattice}
      \index{semilattice!translated --|uu}%
      \index{translated semilattice|uu}%
      if it spans $V$, is discrete and satisfies
      \begin{equation*}
        S+2S\subseteq S.
      \end{equation*}
      A~translated semilattice is called a 
      \emph{semilattice}\index{semilattice|uu} 
      if it contains zero.
    \end{defn}
    
    Because of its significance for the structure of extended affine root
    systems we repeat their classification obtained in \cite{EALA_AMS}.
    
    \begin{constr}\label{constr:eARS}
      Suppose that $\dot R$ is an irreducible finite root system of type
      $X$ in a finite dimensional real vector space $\dot\calV$ with
      positive definite symmetric bilinear form $(\cdot,\cdot)$. We
      decompose the set $\dot R^\times$ of nonzero elements of $\dot R$
      according to length as 
      \begin{equation}\label{eq:shLgEx}
        \dot R^\times
        =\dot R_\text{sh}
        \dot\cup
        \dot R_\text{lg}
        \dot\cup
        \dot R_\text{ex}.
      \end{equation}
      Let $\calV^0$ be a finite dimensional real vector space, set
      $\calV:=\dot\calV\oplus\calV^0$, and extend $(\cdot,\cdot)$ to 
      $\calV$ in such a way that $(\calV,\calV^0)=\{0\}$.
      \begin{enumerate}
  \item[(a)]
        (\emph{The simply laced construction})~
        \index{simply laced construction|uu}%
        \index{construction!simply laced --|uu}%
        Suppose that $X$ is simply laced, i.e. 
        \begin{equation*}
          X=
          A_\ell(\ell\ge1),~
          D_\ell(\ell\ge4),~
          E_6,~
          E_7,~\text{or}~
          E_8.
        \end{equation*}
        Suppose that $S$ is a semilattice in $\calV^0$.
        If $X\neq A_1$ suppose further that $S$ is a lattice in
        $\calV^0$. Put
        \begin{equation*}
          R=R(X,S)
          :=
          (S+S)
          \cup
          \Big(\bigcup_{\dot\alpha\in\dot R^\times}
          (\dot\alpha+S)\Big).
        \end{equation*}
  \item[(b)]
        (\emph{The reduced non-simply laced construction})~
        \index{reduced non-simply laced construction|uu}%
        \index{construction!reduced non-simply laced --|uu}%
        Suppose that $X$ is reduced and non-simply laced, i.e. 
        \begin{equation*}
          X=
          B_\ell(\ell\ge2),~
          C_\ell(\ell\ge3),~
          F_4,~\text{or}~
          G_2.
        \end{equation*}
        Suppose that $S$ and $L$ are semilattices in
        $\calV^0$ such that
        \begin{equation*}
          L+kS\subseteq L
          \text{~~and~~}
          S+L\subseteq S,
        \end{equation*}
        where $k$ is defined by
        \begin{equation}\label{eq:kCases}
          k:=
          \begin{cases}
            2&
            \text{if }R\text{ has type }
            B_\ell(\ell\ge2),~
            C_\ell(\ell\ge3),\\
            &~~~~~~~
            F_4,~
            \text{or}~
            BC_\ell(\ell\ge2),~\\
            3&
            \text{if }R\text{ has type }
            G_2.
          \end{cases}
        \end{equation}
        Further, 
        if $X=B_\ell(\ell\ge3)$ suppose that $L$ is a lattice,
        if $X=C_\ell(\ell\ge3)$ suppose that $S$ is a lattice,
        and if $X=F_4$ or $G_2$ suppose that both $S$ and
        $L$ are lattices. Put
        \begin{equation*}
          R=R(X,S,L)
          :=
          (S+S)
          \cup
          \Big(\bigcup_{\dot\alpha\in\dot R_\text{sh}}
          (\dot\alpha+S)\Big)
          \cup
          \Big(\bigcup_{\dot\alpha\in\dot R_\text{lg}}
          (\dot\alpha+L)\Big).
        \end{equation*}
  \item[(c)]
        (\emph{The $BC_\ell$ construction, $\ell\ge 2$})~
        \index{BC@$BC_\ell$ construction|uu}%
        \index{construction!$BC_\ell$ --|uu}%
        Suppose that $X=BC_\ell(\ell\ge2)$.
        Suppose that $S$ and $L$ are semilattices in
        $\calV^0$ and $E$ is a translated semilattice in
        $\calV^0$ such that $E\cap2S=\emptyset$ and
        \begin{eqnarray*}
          &&
          L+2S\subseteq L,~~
          S+L\subseteq S,\\
          &&      
          E+2L\subseteq E,~~
          L+E\subseteq L.
        \end{eqnarray*}
        If $\ell\ge3$, suppose further that $L$ is a lattice.
        Put
        \begin{eqnarray*}
          R=R(BC_\ell,S,L)
          &:=&
          (S+S)
          \cup
          \Big(\bigcup_{\dot\alpha\in\dot R_\text{sh}}
          (\dot\alpha+S)\Big)\\
          &&
          \cup
          \Big(\bigcup_{\dot\alpha\in\dot R_\text{lg}}
          (\dot\alpha+L)\Big)
          \cup
          \Big(\bigcup_{\dot\alpha\in\dot R_\text{ex}}
          (\dot\alpha+E)\Big).
        \end{eqnarray*}
  \item[(d)]
        (\emph{The $BC_1$ construction, $\ell\ge 2$})~
        \index{BC@$BC_1$ construction|uu}%
        \index{construction!$BC_1$ --|uu}%
        Suppose that $X=BC_1$.
        Suppose that $S$ is a semilattices in
        $\calV^0$ and $E$ is a translated semilattice in
        $\calV^0$ such that $E\cap2S=\emptyset$ and
        \begin{equation*}
          E+4S\subseteq E~~\text{and}~~
          S+E\subseteq S.
        \end{equation*}
        Put
        \begin{eqnarray*}
          R=R(BC_1,S,E)
          &:=&
          (S+S)
          \cup
          \Big(\bigcup_{\dot\alpha\in\dot R_\text{sh}}
          (\dot\alpha+S)\Big)\\
          &&
          \cup
          \Big(\bigcup_{\dot\alpha\in\dot R_\text{ex}}
          (\dot\alpha+E)\Big).
        \end{eqnarray*}
      \end{enumerate}
    \end{constr}
 
    \begin{rem}
      Note that for arbitrary subsets $A$ and $B$ of a vector space,
      the condition $A+B\subseteq A$ is equivalent to 
      $A+\langle B\rangle\subseteq B$ where $\langle B\rangle$ stands
      for the additive subgroup of the vector space generated by
      $B$. This means, for instance that $S+2S\subseteq S$ in the
      definition of a semilattice is equivalent to 
      $S+2\langle S\rangle\subseteq S$. Moreover, the condition
      $L+2S\subseteq L$ in part (b) of the previous construction is
      equivalent to $L+2\langle S\rangle\subseteq L$, etc.
    \end{rem}
    
    \begin{thm}\label{thm:classEARS}
      Let $X$ be one of the types for a finite root system.
      Starting from a finite root system $\dot R$ of type $X$ and up to
      three semilattices or translated semilattices (as indicated in the
      construction), Construction~\ref{constr:eARS} produces an extended
      affine root system of type $X$. Conversely, any extended affine root
      system of type $X$ is isomorphic to a root system obtained from the
      part of Construction~\ref{constr:eARS} corresponding to type $X$.
    \end{thm}

    We need to do some preparations for the definition of the extended
    affine Weyl group of an EARS $R$. Set
    \begin{equation}\label{eq:vV}
      V:=\calV^0\oplus\dot\calV\oplus(\calV^0)^*
    \end{equation}
    and define a bilinear form $\langle\cdot,\cdot\rangle$
    on $V$ in the following way:
    \begin{itemize}
\item 
      $\langle\cdot,\cdot\rangle$ extends $(\cdot,\cdot)$ on
      $\calV=\calV^0\oplus\dot\calV$, 
\item
      $\big\langle\dot\calV,(\calV^0)^*\big\rangle=\{0\}$,~~
      $\big\langle(\calV^0)^*,(\calV^0)^*\big\rangle=\{0\}$,
\item
      $\langle\cdot,\cdot\rangle$ is the natural pairing on
      $\calV^0\times(\calV^0)^*$.
    \end{itemize}
    In other words we have
    \begin{equation*}
      \big\langle(v_1,v_2,v_3),(u_1,u_2,u_3)\big\rangle
      =u_3(v_1)+(v_2,u_2)+v_3(u_1)
    \end{equation*}
    with respect to the decomposition (\ref{eq:vV}).
    Now we have the vector space $V$ with a non-degenerate symmetric
    bilinear form $\langle\cdot,\cdot\rangle$ and a
    totally isotropic subspace $\calV^0$. 

    For ever element $\alpha\in R^\times$ we define a reflection via
    \begin{equation*}
      r_\alpha:~V\to V,
      r_\alpha(v)
      =v-2\frac{\langle v,\alpha\rangle}%
      {\langle\alpha,\alpha\rangle}\alpha
      =v-\langle v,\alpha^\vee\rangle\alpha
      \text{~~~with~~~} 
      \alpha^\vee
      :=
      \frac2{\langle\alpha,\alpha\rangle}\alpha. 
    \end{equation*}  
    \begin{defn}
      The 
      \emph{extended affine Weyl group (EAWeG)}
      \index{group!extended affine Weyl --|uu}%
      \index{extended affine!-- Weyl group|uu}%
      $\calW$
      is the subgroup of $\Aut(V)$ generated by the reflections
      $\{r_\alpha:\alpha\in R^\times\}$. 
    \end{defn}

    In the following result we study the orbits under an EAWeG $\calW$
    of certain points in the span of the root system. 
    This result is need in
    the proof of one of our main theorems.

    \begin{prop}\label{prop:orbits}
      Let $\alpha\in\calV$ and $\dot\alpha\in\dot R$ with
      $\alpha-\dot\alpha\in\calV^0$. We set
      \begin{equation*}
        T
        :=
        \big(\dot\alpha,\dot R^\vee_\sh \big)\langle S\rangle
        +\big(\dot\alpha,\dot R^\vee_\lng\big)\langle L\rangle
        +\big(\dot\alpha,\dot R^\vee_\ex \big)\langle E\rangle.
      \end{equation*}
      Then
      \begin{equation}\label{eq:orbit}
        \calW.\alpha
        =
        \alpha-\dot\alpha+\dot\calW\dot\alpha
        +T
      \end{equation}
    \end{prop}
    \begin{proof}
      To prove that a set is an orbit of a group, two things need to
      be proved: The set is invariant under the group and the group
      acts transitively on the set.

      For invariance we suppose that 
      $\alpha-\dot\alpha+\dot w.\dot\alpha+\tau$ is an
      element of the set on the right side of
      (\ref{eq:orbit}), i.e. 
      $\dot w\in\dot\calW$ and 
      $\tau\in T$. Let $\beta\in R$. Then it can be written as
      $\beta=\dot\beta+\sigma$ with
      $\dot\beta\in\dot R$ and $\sigma\in\calV^0$.

      \begin{eqnarray*}
        r_{\dot\beta+\sigma}
        (\overbrace{\alpha-\dot\alpha}^{\in\calV^0}+\dot w.\dot\alpha+\tau)
        &=&
        \alpha-\dot\alpha+\dot w.\dot\alpha+\tau
        -(\dot w.\dot\alpha,\dot\beta^\vee)
        (\dot\beta+\sigma)\\
        &=&
        \alpha-\dot\alpha
        +\underbrace{r_{\dot\beta}\dot w.\dot\alpha}%
        _{\in\dot\calW.\dot\alpha}
        +\underbrace{\tau}_{\in T} 
        -\underbrace{(\dot\alpha,\dot w^{-1}.\dot\beta^\vee)\sigma}_{\in T}
      \end{eqnarray*}
      Thus, since $\calW$ is generated by the reflections associated to
      the roots, the right hand side of (\ref{eq:orbit}) is invariant
      under $\calW$.

      Now we prove that $\calW$ acts transitively on this set.
      It suffices to show that for
      any $\tau\in T$, any $\dot\beta\in\dot R^\times$ and any
      $\sigma\in S_{\dot\beta}=\{\rho\in\calV^0:\dot\beta+\rho\in R\}$ 
      there is a $w\in\calW$ such that
      \begin{equation*}
        w.(\alpha+\tau)=\alpha+\tau+(\alpha,\beta^\vee)\sigma.
      \end{equation*}
      We have
      \begin{eqnarray*}
        r_{\dot\beta}r_{\dot\beta-\sigma}.(\alpha+\tau)
        &=&
        r_{\dot\beta}\left(
          \alpha+\tau-(\dot\alpha,\dot\beta^\vee)(\dot\beta-\sigma)
        \right)\\
        &=&
        r_{\dot\beta}\left(
          r_{\dot\beta}.\alpha
          +\tau+(\dot\alpha,\dot\beta^\vee)\sigma
        \right)\\
        &=&
        \alpha+\tau+(\dot\alpha,\dot\beta^\vee)\sigma
      \end{eqnarray*}
      This concludes the proof.
    \end{proof}

    \begin{defn}
      Let $R^\times$ be a subset of $\calV$. Then the set
      \begin{equation*}
        \IRC(R^\times):=\big((R^\times-R^\times)\cap\calV^0\big)
        \cup R^\times
      \end{equation*}
      is called the
      \emph{isotropic root closure}.
    \end{defn}

    \begin{lem}\label{lem:iRC}
      If $R$ is an EARS and $R^\times$ is the set of its anisotropic
      roots then $\IRC(R^\times)=R$.
    \end{lem}

    \begin{proof}
      The inclusion $\supseteq$ follows from axiom (R6). For the
      converse inclusion let
      $\sigma\in\big((R^\times-R^\times)\cap\calV\big)$. So there are
      $\alpha$ and $\beta\in R^\times$ with
      $\sigma=\beta-\alpha\in\calV^0$.
      This implies $(\alpha,\beta)=(\alpha,\alpha)>0$. So by (R6) the
      set $\{\beta+n\alpha:n\in\Z\}\cap R$ contains
      $\sigma=\beta-\alpha$, since 
      $d=2\frac{(\alpha,\beta)}{(\alpha,\beta)}+u>u\ge0$. This implies
      $\sigma\in R$, and we are done.
    \end{proof}

    \begin{defn}
      An EARS $R$ is called 
      \emph{minimal}
      if the root set $R^\times$ for the Weyl group $\calW$ 
      is a minimal.
    \end{defn}

    \begin{exmp}\label{exmp:nonMinimal}
      Consider the real vector space $V:=\R^7=\R^3\times\R\times\R^3$ 
      with the bilinear form
      $\langle\cdot,\cdot\rangle$ with matrix
      \begin{equation*}
        \begin{pmatrix}
          0&0&\mathbf{1}\\
          0&1&0\\
          \mathbf{1}&0&0
        \end{pmatrix}
      \end{equation*}
      with respect to the decomposition
      $V:=\R^7=\R^3\times\R\times\R^3$. So $\mathbf{1}$ stands for the
      $3{\times}3$-identity matrix.
      In $V$ we consider the set
      \begin{eqnarray*}
        R^\times
        &:=&\Z^3\times\{-1,1\}\times\{0\}^3.
      \end{eqnarray*}
      This is the set of anisotropic roots of the EARS $\IRC(R^\times)$ 
      of type $A_1$ with nullity 3. 
      Let $\calW$ be the Weyl group of this EARS.

      Set
      \begin{equation*}
        \gamma:=(1,1,1,1,0,0,0)^T\in R^\times.
      \end{equation*}
      Recall the definition of $R^\times_\gamma$ 
      from (\ref{eq:rGamma}). Then 
      \begin{eqnarray*}
        R^\times_\gamma
        &=&
        R^\times\setminus
        \Big(\{-\gamma,\gamma\}
        +\big((2\Z)^3\times\{0\}^4\big)\Big)\\
      &=&
        \{(z_i)_{i=1}^7\in R^\times:z_1z_2z_3~\text{even}\}
      \end{eqnarray*}
      and this is the set of anisotropic roots of an EARS 
      of type $A_1$ with nullity 3.

      We will show that the reflections associated to $R^\times_\gamma$
      generate $\calW$, which means that $R$ as not a minimal
      EARS. This can be done by realizing that $\gamma$ is a so-called
      ghost root for the EARS $\IRC(R^\times_\gamma)$ 
      (see \cite{myThesisBook}~Example~4.4.83). Or it can be done by a
      making
      direct computation.
      For this purpose let the roots
      $\alpha_0,~\alpha_1,\dots,\alpha_6\in R_\gamma$be given by the
      columns of the matrix
      \begin{equation*}
       \begin{pmatrix}
         0	&1	&0	&-1	&0	&-1	&0\\
         -1	&0	&0	&1	&-1	&0	&0\\
         -1	&-1	&1	&0	&0	&0	&0\\
         1&1&1&1&1&1&1\\
         0&0&0&0&0&0&0\\
         0&0&0&0&0&0&0\\
         0&0&0&0&0&0&0
       \end{pmatrix}
      \end{equation*}
      with respect to the standard basis. Then a direct computation shows
      that
      \begin{equation*}
        r_\gamma=
        r_{\alpha_0}
        r_{\alpha_1}
        r_{\alpha_2}
        r_{\alpha_3}
        r_{\alpha_4}
        r_{\alpha_5}
        r_{\alpha_6}.
      \end{equation*}
      This implies that the reflection associated to $R^\times_\gamma$
      generate $\calW$.
    \end{exmp}

    In the proof of the next theorem we will need to handle EARS of
    type $BC_\ell$ with special care. They are best dealt with by
    trimming them into an EARS of reduced type. This motivates the
    following definition.
    
    \begin{defn}
      Let $R$ be an EARS of non-reduced type, i.e. of type
      $BC_\ell$, $\ell\ge1$. We set
      \begin{eqnarray*}
        \trim(R)
        &=&
        \IRC(
        R_\sh\cup R_\lng\cup\frac12 R_\ex)
      \end{eqnarray*}
      and call $\trim(R)$ 
      the \emph{trimmed root system} of $R$.
    \end{defn}

    The following result justifies this terminology.

    \begin{lem}\label{lem:trimmed}
      \begin{enumerate}
  \item The set $\trim(R)$ is an EARS of type $B_\ell$.
  \item The root system $\trim(R)$ has the same Weyl group as $R$.
  \item The root system $\trim(R)$ has the same group $\hat\calW$ (see
        (\ref{eq:preConj})) as $R$. In particular $\calW$ has a
        presentation by conjugation 
        with respect to $R$ if and only if it has a
        presentation by conjugation 
        with respect to $\trim(R)$.
  \item The root system $\trim(R)$ is minimal if and only if $R$ is
        minimal.  
      \end{enumerate}
    \end{lem}

    \begin{proof}
      In order to show that $\trim(R)$ is an EARS it suffices to
      see that $R_\sh\cup R_\lng\cup\frac12 R_\ex$ is the set of
      anisotropic roots of an EARS in view of Lemma~\ref{lem:iRC}.
      If we set $S':=S\cup\frac12E$, we have
      \begin{equation}\label{eq:trR}
        R_\sh\cup R_\lng\cup\frac12 R_\ex
        =
        \big(S'+S'\big)
        \cup
        \left(
          \bigcup_{\dot\alpha\in\dot R_\sh}(\dot\alpha+S')
        \right)
        \cup
        \left(
          \bigcup_{\dot\alpha\in\dot R_\lng}(\dot\alpha+L)
        \right).
      \end{equation}
      A look at Theorem~\ref{thm:classEARS} reveals that it 
      suffices to show that $S'$ is a semi lattice satisfying
      \begin{equation*}
        L+2S'\subseteq L
        ~~~\text{and}~~~
        S'+L\subseteq S'.
      \end{equation*}
      We have
      \begin{eqnarray*}
        S'+2S'
        &=&
        (S\cup\frac12E)+2(S\cup\frac12E)\\
        &=&
        \underbrace{(S+2S)}_{\subseteq S}
        \cup
        \underbrace{(S+E)}_{\subseteq S}
        \cup
        \underbrace{(\frac12E+2S)}_{\subseteq\frac12E}
        \cup
        \underbrace{(\frac12E+E)}_{\subseteq\frac12E}\\
        &\subseteq&
        S\cup\frac12E=S',\\
        L+2S'
        &=&
        L+2(S\cup\frac12E)\\
        &=&
        \underbrace{(L+2S)}_{\subseteq L}
        \cup
        \underbrace{(L+E)}_{\subseteq L}\\
        &\subseteq&
        L
        ~~~\text{and}\\
        S'+L
        &=&
        (S\cup\frac12E)+L\\
        &=&
        \underbrace{(S+L)}_{\subseteq S}
        \cup
        \underbrace{(L+\frac12E)}_{\subseteq\frac12E}\\
        &\subseteq&
        S\cup\frac12E=S'
      \end{eqnarray*}
      From (\ref{eq:trR}) it is clear that $\trim(R)$ is of the type 
      that you obtain by omitting the extra long roots in a $BC_\ell$
      type finite root system. This is type $B_\ell$. 

      Item (ii) follows from the fact that the two Weyl groups have
      the same generating set and (iii) and (iv) 
      follow from the one-to-one
      correspondence of orbits of the Weyl group in $R$ and $\trim(R)$.
    \end{proof}

    \begin{thm}\label{thm:minimal}
      Let $R$ be an EARS and let $\calW$ be its Weyl group. 
      Then $\calW$ has a
      presentation by conjugation with respect to $R$ 
      if and only if $R$ is minimal. 
    \end{thm}

    \begin{proof}
      Let $R$ be an EARS and let $\calW$ be its Weyl group. In view of
      Lemma~\ref{lem:trimmed}, we may assume that $R$ is of reduced
      type.

      In \cite{EAWG}
      a subset $R_\Pi$ of $R$ is introduced which is also an EARS (See
      Proposition~4.12 of \cite{EAWG}). Furthermore it is proven that
      the Weyl group $\calW_\Pi$ of that root system $R_\Pi$ has a
      presentation by 
      conjugation (See Theorem~5.15 of \cite{EAWG}). If we prove that
      \begin{equation}\label{eq:req}
        (\forall w\in\calW)
        ~
        (\forall \alpha\in R^\times)
        ~
        (\exists w'\in\calW_{R_\Pi})
        ~
        w.\alpha=w'.\alpha
      \end{equation}
      then we are done by Proposition~\ref{prop:SubRootPres}.

      Due to Proposition~4.41 and Equation (4.16) in \cite{EAWG} we
      only need to consider EARSs $R$ of type 
      $X=A_1$, 
      $X=B_\ell~(\ell\ge 2)$ and
      $X=C_\ell~(\ell\ge 3)$. According to \cite{EAWG}~(4.11), if we have
      \begin{eqnarray*}
        R&=&
        \big(S+S)
        \cup
        (\dot R_\sh+S)
        \cup
        (\dot R_\lng+L)
        \text{~~~then}\\
        R_\Pi&=&
        \big(S_\Pi+S_\Pi\big)
        \cup
        (\dot R_\sh+S_\Pi)
        \cup
        (\dot R_\lng+L_\Pi)
        \text{~~~with}\\
        &&
        \langle S_\Pi\rangle=\langle S\rangle
        \text{~~~and}~~~
        \langle L_\Pi\rangle=\langle L\rangle
      \end{eqnarray*}
      In view of Proposition~\ref{prop:orbits} this means that the two
      Weyl groups $\calW$ and $\calW_\Pi$ have the same orbits in 
      $R^\times$.
      This entails our requirement in (\ref{eq:req}).
    \end{proof}

    The following characterization of EARSs will be usefull when we
    need to identify certain subsets of EARSs as EARSs.

    \begin{prop}\label{prop:eARSChar}
      Let $R^\times$ be a set of anisotropic vectors in $\calV$ such that
      \begin{enumerate}
  \item $\calW_{R^\times}.R^\times\subseteq R^\times$,
  \item $\overline{R^\times}$ is an irreducible finite root system,
  \item $\langle R^\times\rangle$ is a lattice in $\calV$,
  \item $\alpha\in R^\times$ implies $2\alpha\notin R^\times$.
      \end{enumerate}
      Then $R:=\IRC(R^\times)$ is the EARS having $R^\times$ 
      as its set of anisotropic
      roots. 
    \end{prop}

    \begin{proof}
      We need to verify the axioms (R1) - (R8).
      Axiom (R1) follows from the definition of $\IRC(R^\times)$.
      Due to the fact that
      $r_\gamma.\gamma=-\gamma$ for every $\gamma\in R^\times$ we have
      $-R^\times=R^\times$. By the definition of $\IRC(R^\times)$,
      this implies (R2).
      Item (iii) implies (R3) and (iv) implies (R4). Since $R$ is
      contained in $\langle R^\times\rangle$, it is discrete, so (R5)
      is 
      satisfied. Since $\overline{\dot R^\times}$ is irreducible (R7)
      holds.
      Axiom (R8) follows from the definition of $\IRC(R^\times)$.

      The only axiom remaining to be proved is (R6). 
      So let $\alpha\in R^\times$ and $\beta\in R$. We will look at
      three different cases.

      {\it Case $\overline\alpha$ and $\overline\beta$ linearly
        independent:} 
      Then let $\calW_{\alpha,\beta}$ be the subgroup of $\calW$
      generated by $\alpha$ and $\beta$. This subgroup stabilizes the
      subspace $\calV_{\alpha,\beta}=\spann(\alpha,\beta)$ of $\calV$. Set
      $R_{\alpha,\beta}=\calW_{\alpha,\beta}.\{\alpha,\beta\}$. Since
      \begin{equation*}
        2\frac{(\gamma,\delta)}{(\delta,\delta)}
        =
        2\frac{(\overline\gamma,\overline\delta)}%
        {(\overline\delta,\overline\delta)}
        \in\Z
      \end{equation*}
      for all $\delta$, $\gamma\in R_{\alpha,\beta}$. So the set
      $R_{\alpha,\beta}$ is a finite root system in
      $\calV_{\alpha,\beta}$. In this situation (R6) follows from
      \cite{bou}~Ch.~VI~$n^0$~1.3~Proposition~9. 

      {\it Case $\beta\in R^0$:} 
      Then
      \begin{equation*}
        2\frac{(\alpha,\beta)}{(\alpha,\alpha)}=0
        ~~~\text{and}~~~
        r_\alpha.(\beta+n\alpha)
        =
        \beta+n\alpha-2\frac{(n\alpha,\alpha)}{\alpha,\alpha)}\alpha
        =
        \beta-n\alpha
      \end{equation*}
      for any $n\in\Z$. As a result of this symmetry and due to the
      fact that 
      \begin{equation*}
        \overline{\{\beta+n\alpha:n\in \Z\cap R^\times\}}
        =
        \{n\overline\alpha:n\in \Z\}\cap\overline R^\times
      \end{equation*}
      is a subset of the finite root system $\overline R^\times$ there
      can only be two cases:
      \begin{equation*}
        \{\beta+n\alpha:n\in \Z\}\cap R
        =
        \begin{cases}
          \{\beta-\alpha,\beta,\beta+\alpha\}&\text{or}\\
          \{\beta-2\alpha,\beta-\alpha,\beta,
          \beta+\alpha,\beta+2\alpha\}.
        \end{cases}
      \end{equation*}
      In both cases (R6) is satisfied.

      {\it Case $\beta\in R^\times$ and $\overline\alpha$ and
        $\overline\beta$ linearly dependent:} 
      Since 
      \begin{equation*}
        \{\beta+n\alpha:n\in \Z\}\cap R
        =
        \{\beta+n(-\alpha):n\in \Z\}\cap R
        \text{~~~and~~~}
        2\frac{(-\alpha,\beta)}{(-\alpha,-\alpha)}
        =
        -2\frac{(\alpha,\beta)}{(\alpha,\alpha)},
      \end{equation*}
      We only need to consider the following cases:
      \begin{enumerate}
  \item $\overline\alpha=\overline\beta$,
  \item $2\overline\alpha=\overline\beta$ and
  \item $\overline\alpha=2\overline\beta$.
      \end{enumerate}
      For each of these cases, it suffices to show that there is an
      $n\in\Z$ such that $\beta':=\beta+n\alpha\in R^0$, since 
      \begin{eqnarray*}
        \{\beta+n\alpha:n\in \Z\}\cap R
        &=&
        \{\beta'+n\alpha:n\in \Z\}\cap R\\
        \{\beta-d\alpha,\dots,\beta+u\alpha\},
        &=&
        \{\beta'-(d+n)\alpha,\dots,\beta'+(u-n)\alpha\}
        \text{~~and~~}\\
        (d+n)-(u-n)=d-u+2n
        &=&2\frac{(\alpha,\beta)}{(\alpha,\alpha)}
        +2\frac{(\alpha,n\alpha)}{(\alpha,\alpha)}
        =2\frac{(\alpha,\beta')}{(\alpha,\alpha)}.
      \end{eqnarray*}
      So the proof of this case is reduced to case discussed previously.

      Case~(i):  
      Then $\beta-\alpha\in R^0$ by the definition of
      $\IRC(R^\times)$.

      Case~(ii):
      Then
      \begin{equation*}
        r_\beta(-\alpha)
        =-\alpha-2\frac{(-\alpha,\beta)}{(\beta,\beta)}\beta
        =-\alpha-2\frac{(-\frac12\overline\beta,\overline\beta)}%
        {(\overline\beta,\overline\beta)}\beta
        =\beta-\alpha\in R.
      \end{equation*}
      Then by the definition of $\IRC(R^\times)$, we have
      $\beta-2\alpha\in R^0$.

      Case~(iii):
      We have
      \begin{equation*}
        r_\alpha(\beta)
        =\beta-2\frac{(\beta,\alpha)}{(\alpha,\alpha)}\alpha
        =\beta-\alpha.
      \end{equation*}
      So
      \begin{equation*}
        \{\beta+n\alpha:n\in \Z\}\cap R
        \supseteq\{\beta,\beta-\alpha\}.
      \end{equation*}
      We see that we must have
      equality, since any element in the left hand side and not
      in the right hand side would have at least three times the lenth
      of $\alpha$.
    \end{proof}

    \begin{lem}\label{lem:subEARS}
      Let $R$ be an EARS with EAWeG $\calW$. Suppose there is a root
      $\beta\in R^\times$ such that the reflections associated to
      \begin{equation*}
        \widetilde{R^\times}:=R^\times\setminus\calW.\beta
      \end{equation*}
      generate $\calW$. Then there is an EARS $\tilde R$, whose
      anisotropic roots are $\widetilde{R^\times}$. When passing from
      $R$ to $\tilde R$, only the following changes in types are
      possible:
      \begin{equation*}
        BC_1\to A_1,~~
        BC_2\to B_2,~~
        BC_\ell\to B_\ell,~~
        BC_\ell\to C_\ell,
        \text{~~where $\ell\ge3$}.
      \end{equation*}

    \end{lem}

    \begin{proof}
      We use Proposition~\ref{prop:eARSChar} for $\tilde{R^\times}$. 
      So we need to verify
      conditions (i) - (iv) of that proposition. Conditions (i),
      and (iv) are evident so first we focus on (iii). 
      Discreteness of $\langle\tilde{R^\times}\rangle$ follows from
      the fact that $\langle R^times\rangle$ is discret. (This csn be
      seen, for instance from Theorem\ref{thm:classEARS}.) 
      Now we suppose that $\spann(\tilde{R^\times})$ is a real
      subset of $U$ of $\calV$. This means that $U$ 
      is invariant under $\calW_{\tilde{R^\times}}$ and thus also
      under $\calW_{R^\times}$. That is a contradiction to axiom (R7)
      for the EARS $R$.

      Now we verify (ii). The group $\overline{\calW_{R^\times}}$ 
      obtained by reducing the elements of the
      group $\calW_{R^\times}$ 
      to automorphisms of $\Aut(\calV/\calV^0)$ is the same
      as the subgroup $\calW_{\overline{R^\times}}$ of
      $\Aut(\calV/\calV^0)$ generated by the reflections associated to the
      roots in $\overline{R^\times}$. So by hypothesis we have 
      \begin{equation*}
        \calW_{\overline{R^\times}}=\calW_{\overline{\tilde{R^\times}}}
      \end{equation*}
      
      By Lemma~\ref{lem:finRootSub}, the set
      $\overline{\tilde{R^\times}}$ is an irreducible finite
      crystallographic root
      system. The restrictions about the type of $\tilde R$ follow
      from this lemma, as well.
    \end{proof}

    \begin{thm}\label{thm:subEARS}
      Every EARS $R$ contains a minimal one $R'$. When passing from
      $R$ to $R'$ only the following changes in type are possible: 
      \begin{equation*}
        BC_1\to A_1,~~
        BC_2\to B_2,~~
        BC_\ell\to B_\ell,~~
        BC_\ell\to C_\ell,
        \text{~~where $\ell\ge3$}.
      \end{equation*}
    \end{thm}

    \begin{proof}
      Due to Proposition~\ref{prop:orbits} the Weyl group $\calW$ of
      $R$ has only a finite number of orbits in $R$. So
      Lemma~\ref{lem:subEARS} can be applied repeatedly until a
      minimal EARS $R'\subseteq R$ is obtained. The only possible
      changes in type are those indicated in Lemma~\ref{lem:subEARS}.  
    \end{proof}

    It is tempting to think that an EARS contained in a minimal EARS
    is minimal. But this is not necessarily the case, as the following
    example reveals:

    \begin{exmp}
      Let $R$ be a minimal EARS of type $A_1$ and nullity 3. (Such an
      EARS exists by Theorem~\ref{thm:classEARS} and
      Theorem~\ref{thm:subEARS}. If $R$ is given by
      \begin{equation*}
        R=(S+S)\cup(\dot R_\sh+S)
      \end{equation*}
      then the set 
      \begin{equation*}
        R=( S'+S')\cup(\dot R_\sh+S') 
        \text{~~~with~~~}
        S'=2\langle S\rangle
      \end{equation*}
      is an EARS contained in $R$ by
      Theorem~\ref{thm:classEARS}. Since this EARS is isomorphic to
      the EARS $R$ discussed in Example~\ref{exmp:nonMinimal}, it is
      not minimal.
    \end{exmp}

    The statements of Theorem~\ref{thm:minimal} and
    Theorem~\ref{thm:subEARS} taken together imply:

    \begin{thm}
      Every EAWeG $\calW$ has the presentation by conjugation with
      respect to some EARS.
    \end{thm}

\newcommand{\btxmarchlong}{March}
\newcommand{\manInPrep}{Manuscript in preparation}

\newcommand{\etalchar}[1]{$^{#1}$}

\end{document}